\documentclass[reqno]{amsart}
\usepackage{epsfig,graphicx}

\usepackage{amssymb}      
\usepackage{amsmath}      

\newcommand{\ben}{\begin{enumerate}}
\newcommand{\een}{\end{enumerate}}
\newcommand{\be}{\begin{equation}}
\newcommand{\ee}{\end{equation}}
\newcommand{\bea}{\begin{eqnarray}}
\newcommand{\eea}{\end{eqnarray}}
\newcommand{\bc}{\begin{center}}
\newcommand{\ec}{\end{center}}

\newtheorem{thm}{Theorem}[section]
\newtheorem{cor}[thm]{Corollary}
\newtheorem{lem}[thm]{Lemma}
\newtheorem{prop}[thm]{Proposition}
\newtheorem{conj}[thm]{Conjecture}

\theoremstyle{definition}
\newtheorem{defn}[thm]{Definition}

\theoremstyle{remark}
\newtheorem{rem}[thm]{\rm\bfseries{Remark}}

\newcommand{\Z}{{\mathbb{Z}}}

\begin{document}   
     
\setcounter{page}{1}

\title[Power subgroups]{On power subgroups of mapping class groups}      
 \author[L.Funar]{Louis Funar}     
\address{Institut Fourier BP 74, UMR 5582, 
University of Grenoble I, 38402 Saint-Martin-d'H\`eres cedex, France}      
\email{louis.funar@ujf-grenoble.fr}

\begin{abstract}
In the first part of this paper we 
prove that the mapping class subgroups generated by the 
$D$-th powers of Dehn twists (with $D\geq 2$) along a sparse 
collection of simple closed curves on an orientable 
surface are right angled Artin groups. 
The second part is devoted to power quotients,  
i.e. quotients by the normal subgroup generated 
by the $D$-th powers of all elements of the mapping class groups. 
We show  first that for infinitely 
many $D$ the power quotient groups are non-trivial. 
On the other hand, if $4g+2$ does not divide 
$D$ then the associated power quotient of the mapping class group of 
the genus $g\geq 3$ closed surface is trivial. 
Eventually, an elementary argument shows that 
in genus 2 there are infinitely many power quotients 
which are infinite torsion groups.

\vspace{0.1cm}

\end{abstract}

\maketitle 

\section{Introduction and statements}

The aim of this paper is to give a sample of results concerning 
power subgroups of mapping class groups. 
We denote by $M(S)$ the {\em mapping class group}  
of the orientable  surface $S$, 
namely the group of isotopy classes of orientation-preserving 
homeomorphisms that fix point-wise the boundary components. 
Set $\Sigma_{g,k}^r$ for the 
orientable surface of genus $g$ with $k$ boundary components and $r$ 
punctures. We will omit the indices $k$ and $r$ in $\Sigma_{g,k}^r$ 
when they are zero.

\begin{defn}
Let $A$ be a collection of (isotopy classes of) 
simple closed curves on the surface $S$. 
We denote by $M(S)(A;D)$ the subgroup generated by 
$D$-th powers of Dehn twists  
along curves in $A$. 
\end{defn}

When $A$ is a set $SCC(S)$ 
of representatives for  
all {\em simple closed curves up to homotopy} 
on the surface $S$ the group $M(S)(SCC(S);D)$ will be denoted  by $M(S)[D]$. 
We will omit the indices $k$ and $r$ in $M(\Sigma_{g,k}^r)[D]$ and 
$M(\Sigma_{g,k}^r)(A;D)$ when they are zero.  
For simplicity, when we do not need to 
specify the surface $\Sigma_{g,k}^r$ we will use the 
notation $M_{g,k}^r$ for $M(\Sigma_{g,k}^r)$ and respectively 
$M_{g,k}^r[D]$ for $M(\Sigma_{g,k}^r)[D]$, with the same convention concerning 
the indices $k$ and $r$, which we omit when they are zero.

Observe that $M_g[D]$ is a normal subgroup of $M_g$, whose definition 
is similar to that of the congruence subgroups of the symplectic groups. 
In fact, let $T_a$ denote the Dehn twist along the simple closed curve $a$. 
Then for every $h\in M_g$ we have $h T_a^Dh^{-1}= T_{h(a)}^D \in M_g[D]$. 
As $M_g[D]$ is generated by the $T_a^D$, for $a$ running over the set of 
all  simple closed curves, it follows that $M_g[D]$ is a normal subgroup. 

\vspace{0.2cm}\noindent
The first results on $M_g[D]$ were obtained by Humphries (\cite{Hum})
who proved that $M_g/M_g[2]$, for each $g\geq 1$, 
$M_2/M_2[3]$ and $M_3/M_3[3]$ are finite, while 
$M_2/M_2[D]$ is infinite when $D\geq 4$.

\vspace{0.2cm}\noindent
On the other hand, using quantum topology techniques we proved in 
\cite{F} that the groups $M_g[D]$ are of infinite index in $M_g$, 
if $g\geq 3$, and $D\not\in\{1,2,3,4,6,8,12\}$. 

\vspace{0.2cm}\noindent
Mapping class groups  have interesting actions on various 
moduli spaces, for instance on spaces of $SU(2)$ 
representations of surface groups. It is known (see \cite{Go}) that 
the whole mapping class group acts ergodically. Actually the same proof 
extends trivially to show that $M_g[D]$ still 
acts ergodically. This yields the first examples of infinite index subgroups
of the mapping class group acting ergodically.

Methods from quantum topology also show that: 
\[ \cap_{D\in \mathcal D} M_g[D]=1\]
if $g\geq 2$ and $\mathcal D$ is  any infinite set of positive integers. 
In fact, the kernel of the $SO(3)$ quantum representation  
of level $k$ of $M_g$ contains $M_g[k]$. Then the 
asymptotic faithfulness theorem 
from \cite{A,FWW} yields the claim.

However, these results seem to exhaust our present knowledge 
about the groups $M_g[D]$. It is not known, for instance, 
whether the following  holds:
\begin{conj}
The group $H_1(M_g[D])$ is not finitely generated if  
$D\geq 3$, $g\geq 4$ or $D\geq 4$, $g\in\{2,3\}$. 
\end{conj}

\vspace{0.2cm}\noindent
If true, this would imply that $M_g/M_g[D]$ is infinite for the above 
values of $D$ and $g$. 

\begin{rem}
The groups $M_g[2]$ have finite index in $M_g$ 
(see \cite{Hum}) and hence are finitely generated.
However the quantum representations at 
$4$-th roots of unity (see \cite{Mas,Wr1}) and 
$6$-th roots of unity (see \cite{Wr2}) have finite image. Thus the 
quantum method used for large $D$ cannot decide whether   
$M_g[4]$ and $M_g[6]$ have finite index or not. 
It is likely that $M_g[D]$ is of infinite index for every $D\geq 4$ and 
$g\geq 3$. Notice also that a similar problem for pure braid groups 
was considered in \cite{Hum2}.  
\end{rem}

A question of Ivanov (see \cite{I3}, Question 12) is particularly relevant 
for the structure of the group $M_g[D]$ by asking about the possible relations 
between powers of Dehn twists. We formulate it here as a 
conjecture, under a slight restriction on $D$: 
 
\begin{conj}\label{braid}
The group  $M_g[D]$ (for $D\geq 3$, $g\geq 4$ or $D\geq 4$, $g\in\{2,3\}$) 
has the following presentation:  
\begin{enumerate}
\item Generators $Z_a$ (standing for $T_a^D$), 
where $a$ belongs to the (infinite) set $SCC(\Sigma_g)$ 
of simple closed curves on the surface; 
\item Relations of conjugacy type: 
\[ Z_{T_a^D(b)}=Z_a Z_b Z_a^{-1}\]
for each pair $a, b\in SCC(\Sigma_g)$. 
\end{enumerate}
\end{conj}

\noindent We denote by ${\mathcal A}_{\Gamma}$ the right angled Artin 
group associated to the graph $\Gamma$, which is defined by the 
following presentation: 
\begin{enumerate}
\item Generators  $Z_a$, where $a$ belongs to the set of vertices 
of $\Gamma$; 
\item Relations 
\[ Z_{a}Z_b=Z_bZ_a, \mbox{ if } a \mbox{ and } b \mbox{ are connected   
by an edge in } \Gamma\]
\end{enumerate}

A related (but much weaker) Conjecture is as follows:

\begin{conj}\label{artin}
Let $C_{\Sigma,D}\subset SCC(\Sigma)$ be a set of representatives 
of  the orbits set $SCC(\Sigma)/M(\Sigma)[D]$. Consider the 
associated intersection graph $\Gamma(C_{\Sigma,D})$, whose vertex set is 
$C_{\Sigma,D}$ and edges join  vertices corresponding to 
disjoint curves on the surface. Then the homomorphism 
$A_{\Gamma(C_{\Sigma,D})}\to M(\Sigma)[D]$ which sends the generators 
$Z_a$ into the elements $T_a^D$ is an isomorphism on its image 
for $D\geq 3$, $g\geq 4$ or $D\geq 4$, $g\in\{2,3\}$. Here $g$ denotes the 
genus of the surface $\Sigma$. 
\end{conj}

Clay, Leininger and Margalit recently proved in \cite{CLM} that $M(\Sigma)[D]$ is not 
abstractly commensurable with any right angled Artin group. In particular, the homomorphism
$A_{\Gamma(C_{\Sigma,D})}\to M(\Sigma)[D]$ from above is not surjective.

\begin{rem}
According to Ishida (see \cite{Ishida}) the group generated by 
two Dehn twists is either free abelian (if the curves are disjoint or coincide) 
or generate the braid group 
$B_3$ in 3 strands (if the curves intersect in one point) or 
free (if the curves intersect in at least two points). 
In particular the subgroup generated by two $D$-th powers of 
Dehn twists is either free abelian or free, supporting the 
Conjecture \ref{artin}. 
See also \cite{CP} or (\cite{Ham}, Thm. 3.5) for the braid case. 
Relations between multi-twists are also  given in \cite{Marg}. 
\end{rem}

\begin{prop}
The analogues of Conjecture \ref{braid}  
for  $D=2$ and any  closed orientable surface $\Sigma$ of genus $g\geq 3$ 
are false as stated, namely there are additional relations in a 
presentation of $M_g[2]$ with the given generators. 
\end{prop} 
\begin{proof} 
According to Humphries (see \cite{Hum}) the subgroup 
$M_g[2]$ can be identified to  the kernel of 
the homomorphism $M_g\to Sp(2g,\Z/2\Z)$. 
Hain  proved in \cite{Ha} (see also \cite{McC}) that 
any finite index subgroup of $M_g$ 
(for $g\geq 3$) containing the Torelli subgroup (i.e. the subgroup of 
mapping classes acting trivially in homology) 
has trivial first cohomology. This implies that $H^1(M_g[2])=0$, 
which was also proved  by McCarthy in \cite{McC}. 
But the abelianization of the group presented by the relations 
from Conjecture \ref{braid}  
is a free abelian group of rank equal to the cardinal of 
$SCC(\Sigma_g)/M_g[2]$. This contradiction shows that in $M_g[2]$ 
there are additional relations.   
\end{proof}

\begin{rem}
The referee pointed out  explicit relations 
among squares of Dehn twists along nonseparating curves. Choose for instance 
the nonseparating curves $a_1,a_2,\ldots,a_7$ on $\Sigma_3$ such that 
$a_i$ intersects $a_j$ at one point if $j=i+1$ and they are disjoint otherwise.
Then we have the following relation in $M_3$:
\[ (T_{a_1}T_{a_2}T_{a_3} T_{a_4}T_{a_5}T_{a_6}T_{a_7}T_{a_6}T_{a_5}T_{a_4}T_{a_3}T_{a_2}T_{a_1})^2=1\]
Observe further that 
\[ T_{a_1}T_{a_2}T_{a_3} T_{a_4}T_{a_5}T_{a_6}T_{a_7}T_{a_6}T_{a_5}T_{a_4}T_{a_3}T_{a_2}T_{a_1}=\]
\[=(T_{a_7}^2)^{A_6}(T_{a_6}^2)^{A_5}(T_{a_5}^2)^{A_4}(T_{a_4}^2)^{A_3}(T_{a_3}^2)^{A_2}(T_{a_2}^2)^{A_1}(T_{a_1}^2)\]
where we put $A_i=T_{a_1}T_{a_2}\cdots T_{a_i}$, and $x^A=AxA^{-1}$. 
Now, we can express $(T_{a_i}^2)^{A_{i-1}}=T_{A_{i-1}(a_i)}^2$ as squares 
of Dehn twists. We obtain therefore  the following relation
\[(T_{A_6(a_7)}^2T_{A_5(a_6)}^2T_{A_4(a_5)}^2T_{A_3(a_4)}^2T_{A_2(a_3)}^2T_{A_1(a_2)}^2T_{a_1}^2)^2=1\]
which does not follow from those defining a right angled Artin group in these (square Dehn twists) generators.  
\end{rem}

\vspace{0.1cm}\noindent 
Remark also that the analogue of Conjecture \ref{braid} cannot hold  
when $D=1$ either. In fact the abelianization  of $M_g$ 
would be a nontrivial free abelian group, 
contradicting the fact that $M_g$ is perfect when $g\geq 3$ and has 
torsion abelianization otherwise.


\vspace{0.2cm}\noindent
An important step towards a solution to Conjecture \ref{artin} was taken 
in the recent paper \cite{Kob} of Koberda, where it is proved the following.  
For any irredundant (see \cite{Kob} for the definition) 
collection $\{f_1,f_2,\ldots,f_k\}$ of mapping 
classes of homeomorphisms, each one being either a Dehn twist or a 
pseudo-Anosov homeomorphism supported on a single connected 
subsurface, there exists $N_0$ such that 
$\{f_1^N,f_2^N,\ldots, f_k^N\}$ is a right angled generating system 
for a right angled Artin subgroup of the mapping class group, 
for any $N\geq N_0$.

The first result of this paper supports  further evidence for the last two   
conjectures.  
Let $A$ be a finite collection of simple closed curves on  a surface 
$S$ and  denote by $F(A)$  the regular neighborhood  of $A$ in $S$. 
We assume that curves are isotoped so that for each $a, b\in A$ 
the number $i(a,b)$ of intersection points between $a$ and $b$ is minimal. 
We pick up a base point $p$ on the 
surface $S$ and a set of distinct points 
$p_a^0\in a$, for $a\in A$. 

\begin{defn}
 The  collection $A$  of curves on the surface $S$ is 
{\em sparse} if it is finite and for some choice 
of paths $\gamma_{a}$ joining $p$ to $p_a^0$ the free 
subgroup $O(A)\subset \pi_1(F(A),p)$ generated 
by the homotopy classes of based loops 
$\gamma_{a}a\gamma_{a}^{-1}$, $a\in A$,  
embeds into $\pi_1(S,p)$ under the map induced by  the inclusion 
$F(A)\hookrightarrow S$.  The collection $A$ is {\em nontrivial} 
if the group $O(A)$ is nontrivial.  
\end{defn}

\begin{thm}\label{main}
Let $D\geq 2$ and $A$ be a nontrivial sparse collection of curves 
on $\Sigma_{g,d}$, where $d\geq 1$. Then the subgroup 
$M(\Sigma_{g,d})(A;D)$ is a right angled Artin group. 
\end{thm}
\begin{rem}
One can construct sparse sets $A$ by considering free subgroups 
generated by primitive elements in $\Sigma_{g,d}$. 
\end{rem}

\begin{rem}
J.Crisp and L.Paris considered before the question of 
finding presentations of subgroups generated by non-trivial 
powers of the standard generators in  Artin groups. 
They established in \cite{CP} the Tits conjecture, which  claimed 
that these subgroups are right angled Artin groups. 
M.L\"onne proved in \cite{Lo} similar results in the braid group 
setting, by showing that the subgroups generated by  the powers of band 
generators are again right angled Artin groups if the powers are at least 3. 
\end{rem}

\begin{rem}
Recently, M.Kapovich proved in \cite{Ka} (making use of our result above) 
that  all right angled Artin groups associated to finite graphs 
embed into the group of Hamiltonian symplectomorphisms of any 
symplectic manifold.   
\end{rem}

\vspace{0.2cm}\noindent
The second part of this article is concerned with power subgroups and 
quotients. Recall the following: 

\begin{defn}
Set  $X_g[D]$ for the $D$-th power subgroup of $M_g$, namely the subgroup 
generated by powers $h^D$ of  all elements of $h\in M_g$. 
Then $X_g[D]$ is a normal subgroup of $M_g$ whose 
quotient is a torsion group.  
\end{defn}

\begin{rem}
Newman (\cite{New}) proved that the $D$-th power subgroup 
of $PSL(2,\Z)$ (and hence of $SL(2,\Z)$) is of infinite index 
when $D=6m\geq 48000$. This was extended by  
Fine and Spellman (see \cite{FS}) to the $2p$-th power subgroups 
of $\Z/2\Z*\Z/p\Z$ (for odd prime $p$).  
 \end{rem}

A natural question is whether power 
quotients of the mapping class group  could  be  non-trivial, 
or even infinite torsion groups. Our second result gives some answers 
in particular cases: 
 
\begin{thm}\label{nontrivial}
\begin{enumerate}
\item Choose an ordered basis of $H_1(\Sigma_g,\Z)$ and denote by 
$P:M_g\to Sp(2g,\Z)$ the homomorphism which sends a mapping class 
into the matrix describing its action in homology.
Then, for every $g\geq 2$ there exist infinitely many  integers $D$ for which  
$P(X_g[D])$ is a proper  
subgroup of $Sp(2g,\Z)$. In particular $M_g/X_g[D]$ are non-trivial 
torsion groups, for these values of $D$. 
\item If $4g+2$ does not divide $D$ and $g\geq 3$ then $M_g=X_g[D]$.
\end{enumerate}
\end{thm}

The question concerning the existence of infinite torsion quotients 
of $M_g$ (see the question of Ivanov in \cite{I3})  
has an elementary solution for genus $g=2$. 
Using arguments similar to those of Korkmaz in \cite{Ko} 
we show that: 

\begin{thm}\label{tors}
The group $M_2/X_2[360 D]$ is an infinite torsion group (of exponent 
$360 D$) for $D\geq 8000$. 
\end{thm}

{\bf Acknowledgements.}  
We are thankful to Thomas Koberda, 
Mustafa Korkmaz, Michael L\"onne, 
Ivan Marin and Luis Paris  for useful discussions and the referee for 
corrections and suggestions. 
The author was partially supported by 
the ANR 2011 BS 01 020 01 ModGroup.

\section{Subgroups of mapping class groups generated by powers of Dehn twists}

\subsection{Finitely generated subgroups generated by powers in braid groups}
The  analogs of the groups $M(\Sigma_g)(A;D)$ in the case of braid groups 
have been considered  long time ago by Coxeter. 
The braid group $B_n$ in $n$ strands  has the usual presentation, due to Artin:           
\[ B_n= \langle \sigma_1,\sigma_2,...,\sigma_{n-1} \mid 
\sigma_i\sigma_j = \sigma_j\sigma_i,\mbox{ if } \mid i-j\mid>1, \;\;
\sigma_{i+1}\sigma_i\sigma_{i+1}=\sigma_i\sigma_{i+1}\sigma_i, 1\leq i\leq n-2 \rangle\]
It is well-known that  the quotient of $B_n$ by 
the normal subgroup generated by $\sigma_i^2$ is the permutation 
group $S_n$. Consider, after Coxeter (see \cite{cox1}): 
\begin{defn} 
The subgroup $B_n\{D\}$ of $B_n$ is the group   
generated by the powers $\sigma_i^D$ of the {\em standard} 
generators $\sigma_i$. Let also $N(B_n\{D\})$ denote the normal 
closure of $B_n\{D\}$ in $B_n$. 
\end{defn}

Coxeter gave  in \cite{cox1} the list of all those quotients   
$B_n/N(B_n\{D\})$ which are finite, together with their respective description
(see also \cite{cox2,cox3}), as follows: 
\begin{prop}[Coxeter]
The group $N(B_n\{D\})$ is of finite index in $B_n$ if and only if 
$(D-2)(n-2) < 4$. Away from the trivial cases $D=2$ or $n=2$ 
we have another five groups:
\begin{enumerate}
\item n=3 
  \begin{enumerate}
  \item For  $D=3$ the quotient $B_3/N(B_3\{3\})$ is isomorphic to  
            $SL(2, \Z/3\Z)$ and has order 24;
  \item For  $D=4$ the quotient $B_3/N(B_3\{4\})$ is a non-split extension of 
the symmetric group $S_4$ on a set of $4$ elements 
by $\Z/4\Z$ and has order 96;
  \item For  $D=5$ the quotient $B_3/N(B_3\{5\})$  is isomorphic to 
 $GL(2, \Z/5\Z)$ and  has order 600; 
  \end{enumerate}
\item For $n=4, D=3$ the factor group $B_4/N(B_4\{3\})$ 
has order 648 and is the central extension of the
Hessian group $(\Z/3\Z)^2\rtimes PSL(2,\Z/3\Z)$ by $\Z/3\Z$.
\item For  $n=5, D=3$ the factor group $B_5/N(B_5\{3\})$ 
has order 155 520 and  is the central extension of the
simple group of order 25 920 by $\Z/6\Z$.
\end{enumerate}
\end{prop}

\begin{rem}
\begin{enumerate}
\item There is an  analogue of the Conjecture \ref{artin}  for the punctured disk
$\Sigma_{0,1}^n$, where we replace powers of 
Dehn twists by  powers of half-twists (i.e. braids).  
Notice that $N(B_n\{2D\})$ is a subgroup of $M_{0,1}^n[D]$.  
\item 
J.Tits conjectured that $B_n\{D\}$ and more generally 
the subgroups generated by powers of  the standard 
generators in Artin groups are right angled Artin groups.   
The later conjecture was settled in full generality by 
Crisp and Paris \cite{CP}. 
\item 
It seems unknown whether the analogue of the Conjecture \ref{artin} 
for  $N(B_n\{D\})$ holds for $D\geq 3$. Notice that for $D=2$ there exist 
nontrivial relations among squares of band generators (which 
are Dehn twists) according to \cite{Lo}. 
\end{enumerate}
\end{rem}

\subsection{Proof of Theorem \ref{main}}

Consider the regular neighborhood $F(A)$ of $A$ in $\Sigma_{g,d}$, which is  a subsurface of genus $g(A)$ with $k(A)$ 
boundary components. Then 
$g(A)\leq g$, but the number $k(A)$ of 
boundary components of $F(A)$ depends on the geometry of $A$ and 
can be arbitrarily large. We denote by  $i(a,b)$  
the minimal number of intersection points between curves 
in the isotopy classes of $a$ and $b$, respectively. 
We assume that curves in $A$ are isotoped so that for each $a, b\in A$ 
the number of intersection points between $a$ and $b$ equals $i(a,b)$ 
and there are not triple intersections among curves in $A$.

We will adapt the proof of the Tits conjecture given in \cite{CP}. 
In the present situation we will be concerned 
with an Artin group  $B(A)$ (to be defined later) 
associated to the collection $A$ and its representation 
into the mapping class group of $F(A)$. 

\vspace{0.2cm}\noindent
We can obtain $F(A)$ as the result of plumbing one  
annulus neighborhood $Ann_a$  for each curve $a$ in $A$.  
In particular the core curves of the annuli are transverse to each other. 
Pick-up one base point $p_a^0$  in the boundary of $Ann_a$, for each $a\in A$.
We can suppose that all $p_a^0$ belong to $\partial F(A)$.  
Choose one distinguished boundary component $a^+$ for each annulus $Ann_a$. 
There is no loss of generality in assuming that each $p_a^0$ belongs to 
$a^+$ and a small arc of $a^+$ centered at $p_a^0$ is contained 
within $\partial F(A)$. 

\vspace{0.2cm}\noindent
Give an orientation to every curve $a\in A$ and a total order $<$ on $A$. 

\vspace{0.2cm}\noindent
If we travel along $a^+$ in the direction given by the orientation 
and starting at $p_a^0$ we will meet a number of intersection points 
between $a^+$ and the other curves $b^+$, where $b\in A$. 
We denote them in order $p_a^1,p_a^2,\ldots, p_a^{d(a)}$.
Denote then by $S=\cup_{a\in A}\cup_{0\leq j\leq d(a)}\{p_a^j\}$ the set of 
all these points. It is clear that $S\subset \partial F(A)$.

\vspace{0.2cm}\noindent
The groupoid $\pi_1(F(A), S)$ is the fundamental groupoid 
of $F(A)$ based at the points of $S$. 
Since $F(A)$ has boundary it follows 
that $\pi_1(F(A), S)$  is a free groupoid (see \cite{DV}, p.7). 

\vspace{0.2cm}\noindent
Furthermore  the mapping class group $M(F(A))$ acts by automorphisms on the 
fundamental groupoid $\pi_1(F(A), S)$, because 
$S\subset \partial F(A)$ and elements of $M(F(A))$ are classes of 
homeomorphisms fixing point-wise the boundary.

\vspace{0.2cm}\noindent
Consider the following elements of $\pi_1(F(A), S)$: 
\begin{enumerate}
\item  For every $s\in A$ the {\em elementary 
loop} $\alpha_s$ is $s^+$ based at $p_s^0$, with its orientation. 
Thus $\alpha_s$  is parallel to the central curve $s$ in the annulus  $Ann_s$. 
\item For every $s\in A$ and $i\in\{0,1,\ldots, d(s)-1\}$ 
consider the arc $p_s^ip_s^{i+1}$ of $s^+$ which 
joins $p_s^i$ to $p_s^{i+1}$. We call them {\em admissible arcs}.  
Observe that the arc $p_s^{d(s)}p_s^0$ is not admissible. 
\end{enumerate}

\vspace{0.2cm}\noindent
The {\em dual intersection graph} of $A$ is defined as the graph whose vertices 
are the elements of $A$ and two vertices are connected by an edge if 
the corresponding curves intersect. 
Assume henceforth that the dual intersection graph of $A$ is connected. 
Then admissible arcs and elementary loops generate the groupoid 
${\mathbb F}=\pi_1(F(A), S)$. 

\vspace{0.2cm}\noindent
Let then $\Gamma_A$ be the subgroup of $M(F(A))$  
generated by the Dehn twists $T_a$, for all $a\in A$. 

\vspace{0.2cm}\noindent
Set $\mathbb B$ for the sub-groupoid of $\mathbb F$ generated by the 
admissible arcs. 

\vspace{0.2cm}\noindent
We will need some terminology  and facts from \cite{CP}. Any element of 
$\mathbb F$ can uniquely be written in the reduced form: 
\[ w=\mu_0 \alpha_{s_1}^{k_1}\mu_1\cdots \alpha_{s_m}^{k_m}\mu_m\]
where $\mu_i\in \mathbb B$, $\mu_i$ is non-trivial if $i\neq 0,m$ and 
$k_i\neq 0$.

We say that $w$ has a {\em square} in $\alpha_s$ if for some $j$ we have  
$s_j=s$ and $|k_j|\geq 2$, and is {\em without squares in} $\alpha_s$,  
otherwise. Moreover $w$ is of type $(\mu, \alpha_t^p)$ if 
its reduced form is 
\[ w=\mu_0 \alpha_{t}^{k_1p}\mu_1\cdots \alpha_{t}^{k_mp}\mu_m, \,\, 
k_j\in \Z\setminus\{0\}, 
\mbox{ and } \mu=\mu_0\mu_1\cdots \mu_m\]
By language abuse we will speak about  
$T_a(w)$, where $w$ is a word in $\mathbb F$, using the action of $\Gamma_A$ 
by automorphisms on $\mathbb F$.

\begin{lem}\label{cal}
Let $s\in A$ and $m\in \Z\setminus\{0\}$. 
\begin{enumerate}
\item If $\mu\in \mathbb B$ then $T_s^m(\mu)$ is of type $(\mu,\alpha_s^m)$. 
\item Let $t\in A$. If $s=t$ or $i(s,t)=0$ 
then $T_s^m(\alpha_t)=\alpha_t$. 
\item If $i(s,t)\neq 0$ then 
$T_t^m(\alpha_s)$ is $u\alpha_s$, where $u$ is an element of 
type $(1,\alpha_t^m)$. 
Thus, if $|m|\geq 2$ and $i(s,t)\neq 0$ 
then $T_t^m(\alpha_s)$ has a square in $\alpha_t$. 
\end{enumerate}
\end{lem}
\begin{proof}
If $s^+,t^+$ intersect at $p$ we define $\varepsilon(s,t; p)\in\{-1,1\}$ 
as follows. Assume that we travel along $s^+$ to meet $p$. At $p$ we use the 
global orientation of the surface for turning right along $t^+$ 
and continue travelling this way. If the direction along $t^+$ 
is the orientation of $t^+$ then we set $\varepsilon(s,t; p)=1$ and 
otherwise  $\varepsilon(s,t; p)=-1$. 

\vspace{0.2cm}\noindent
Next, we will identify
canonically $\pi_1(F(A),S)$ with $\pi_1(F(A), S')$ where $S'$ 
is a copy of $S$, each point $p_a^j$ being slightly moved  
in the positive direction along the arc $a^+$ to a point $\tilde{p}_a^j$. 

\vspace{0.2cm}\noindent
Denote by $\alpha_t(p_t^j)$ the element 
$\tilde{p}_t^0\tilde{p}_t^j \, \alpha_t \, \tilde{p}_t^j \tilde{p}_t^0$, where $\tilde{p}_t^0\tilde{p}_t^j$ is the unique 
arc of $\alpha_t$ joining $\tilde{p}_t^0$ to $\tilde{p}_t^j$  and consisting only 
of admissible sub-arcs. Also $\alpha_t$ is considered as 
the loop $t^+$  based at $\tilde{p}_t^0$. 
Then by  direct computation we find: 
\[ T_{\alpha_t}^m(\tilde{p}_s^i\tilde{p}_s^{i+1})=
\left\{ \begin{array}{ll}
\tilde{p}_s^i\tilde{p}_s^{i+1}, & \mbox{ if } p_s^ip_s^{i+1}\cap \alpha_t =\emptyset, \mbox{ or } 
s=t \\
\tilde{p}_s^i\tilde{p}_s^{i+1}\, \alpha_t^{m\varepsilon(s,t;p_s^{i+1})}(p_s^{i+1}) & 
 \mbox{ if } p_s^{i+1}\in t^+\\
\tilde{p}_s^i\tilde{p}_s^{i+1} & 
\mbox{ if } p_s^{i+1}\not\in t^+ \\
\end{array}\right. 
\]

\vspace{0.2cm}\noindent
Notice that when  the start-point $p_s^i$ belongs to $t^+$ the action 
is trivial since the base-point $p_s^i$ is slightly pushed out of $t^+$ 
in $S'$. 

\vspace{0.2cm}\noindent
Let now $s,t\in A$ be two curves with $i(s,t)\neq 0$.  
Suppose now that starting from $p_s^0$ and  
traveling along $s^+$ we meet the circle $t^+$ at the points 
$p_s^{j_1},p_s^{j_2},\ldots, p_s^{j_r}$, $r>0$. 
By direct inspection we find that  
\[  T_{\alpha_t}^m(\alpha_s)=
\tilde{p}_s^0\tilde{p}_s^{j_1}\, \alpha_t^{m\varepsilon(s,t;p_s^{j_1})}(p_s^{j_1}) \, 
\tilde{p}_s^{j_1}\tilde{p}_s^{j_2}\alpha_t^{m\varepsilon(s,t;p_s^{j_2})}(p_s^{j_2}) \cdots 
 \alpha_t^{m\varepsilon(s,t;p_s^{j_{r-1}})}(p_s^{j_{r-1}})(\tilde{p}_s^{0}\tilde{p}_s^{j_r})^{-1}\alpha_s\]
It is immediate that  $T_{\alpha_t}^m(\alpha_s)=u\alpha_s$, where 
$u$ is of type $(1,\alpha_t^m)$. 
\end{proof}

\begin{lem}\label{square}
Let $x\in \mathbb F$, $|m|\geq 2$. If $x$ is without squares in $\alpha_t$ and 
$T_{\alpha_t}^m(x)$ has a square in $\alpha_s$ then either $s=t$ or else 
$i(s,t)=0$ and $x$ has a square in $\alpha_s$.  
\end{lem}
\begin{proof}
Let $x=\mu_0 \alpha_{s_1}^{k_1}\mu_1\cdots \alpha_{s_r}^{k_r}\mu_r$ 
in reduced form. 
The previous lemma shows that: 
\begin{enumerate}
\item If $s_i=t$ then 
$v_i=T_{\alpha_t}^m(\alpha_{s_i}^{k_i})=\alpha_{s_i}^{k_i}$, where 
$k_i\in\{-1,1\}$, because $x$ is without squares in $\alpha_t$. 
\item If $s_i$ and $t$ are disjoint then 
$T_{\alpha_t}^m(\alpha_{s_i}^{k_i})=\alpha_{s_i}^{k_i}$. 
\item If $i(s_j,t)\neq 0$ then 
\[ T_{\alpha_t}^m(\alpha_{s_j}^{k_j})=
\left\{\begin{array}{ll}
u_j(\alpha_{s_j}u_j)^{k_j-1}\alpha_{s_j} & \mbox{ if } k_j >0\\
\alpha_{s_j}^{-1}(u_j^{-1}\alpha_{s_j}^{-1})^{-k_j-1}u_j^{-1} 
& \mbox{ if } k_j <0 
\end{array}\right.
\] 
where $u_j$ is  a non-constant term of type $(1,\alpha_{t}^m)$. 
\item  $T_{\alpha_t}^m(\mu_j)$ has a reduced form $y_j$ of type 
$(\mu,\alpha_t^m)$, for all $j\geq 0$. 
\end{enumerate}

\vspace{0.2cm}\noindent
Therefore we can write in reduced form 
$T_{\alpha_t}^m(x)= x_0v_1x_1v_2\cdots v_rx_r$ as follows: 
\begin{enumerate}
\item 
If either $s_i=t$ or $s_i$ and $t$ are disjoint 
then $v_i=T_{\alpha_t}^m(\alpha_{s_i}^{k_i})=\alpha_{s_i}^{k_i}$.
\item Assume that $i(s_j,t)\neq 0$. 
\begin{enumerate}
\item If $k_j>0$ then $v_j=(\alpha_{s_j}u_j)^{k_j-1}\alpha_{s_j}$. 
Absorb the extra factor $u_j$ into $x_{j-1}$.  
\item If $k_j<0$ then $v_j=\alpha_{s_j}^{-1}
(u_j^{-1}\alpha_{s_j}^{-1})^{-k_j-1}$. Absorb the extra factor $u_j^{-1}$ 
into $x_{j}$.  
\end{enumerate}
\item Eventually $x_j$ are $T_{\alpha_t}^m(\mu_j)$, possibly corrected by 
the absorption of terms coming from $v_{j}$ or $v_{j+1}$.  
Thus $x_j$ are of reduced form of type $(\mu_j,\alpha_t^m)$. 
\end{enumerate}
In particular, if $T_{\alpha_t}^m(x)$ has a square in $\alpha_s$ 
then either $s=t$ or there exists $j$ such that $s_j=s$ and 
$s$ and $t$ are disjoint. 
\end{proof}

\vspace{0.2cm}\noindent
To each set of curves $A\subset \Sigma_{g,d}$ we  
can associate the Artin group $B(A)$, with the following 
presentation:
\[ B(A)=\langle z_a, a\in A\mid z_az_b=z_bz_a, \mbox { if } 
a\cap b =\emptyset, \, z_az_bz_a=z_bz_az_b,  \mbox{ if } i(a,b)=1\rangle \]
There is a natural homomorphism $\tau:B(A)\to M(F(A))$ 
which sends  $z_a$ into the Dehn twist $T_a$.

\vspace{0.2cm}\noindent
Consider now the right angled Artin group defined by the presentation: 
\[ H(A)=\langle w_a, a\in A\mid w_aw_b=w_bw_a, \mbox { if } 
i(a,b) =0 \rangle \]

\vspace{0.2cm}\noindent
There is a map $\iota:H(A)\to B(A)$ given by $\iota(w_a)=z_a^D$. We 
will suppose that $D\geq 2$ in the sequel.  
The word 
$W=w_{s_l}^{n_l}w_{s_{l-1}}^{n_{l-1}}\cdots w_{s_2}^{n_2}w_{s_1}^{n_1}$ 
is called an {\em M-reduced expression} of the 
element $w\in H(A)$ (obtained by interpreting letters 
as the corresponding generators of $H(A)$) if for any $i<j$ such that 
$s_i=s_j$ there exists $k$ such that $i<k<j$ and $i(s_i,s_k)\neq 0$. 
Then the $M$-reduced expression for $w$ {\em ends} in $s$ if, up 
to change the order of commuting generators, we can arrange that $s_l=s$. 

\vspace{0.2cm}\noindent
Recall now that $\tau(\iota(w))$ is an automorphism of $\mathbb F$, for each 
$w\in H(A)$. We will write simply $w(x)$ or $W(x)$ for $\tau(\iota(w))(x)$, 
where $w\in H(A)$, $x\in \mathbb F$ and $W$ is an $M$-reduced expression for $w$. 

\vspace{0.2cm}\noindent
The following two lemmas are restatements of 
Propositions 9 and 10 from \cite{CP}.

\begin{lem}
Let $W$ be an $M$-reduced expression for $w\in H(A)$, $x\in \mathbb F$ 
and $s\in A$. Suppose that $x$ is without squares in 
$\alpha_t$ for all $t\in A$, 
and that $w(x)$ has a square in $\alpha_s$. Then $W$ ends in $s$. 
\end{lem}
\begin{proof}
We will proceed by induction on the length $l$ of  the $M$-reduced expression 
$W=w_{s_l}^{n_l}w_{s_{l-1}}^{n_{l-1}}\cdots w_{s_2}^{n_2}w_{s_1}^{n_1}$ (see also \cite{CP},p.30). 
When $l=0$, $w$ is identity and  thus $w(x)=x$ cannot have 
squares in $\alpha_s$, under our assumptions. For the induction step 
let now write $W=w_{s_l}^{n_l}W'$, where $l\geq 1$. If $W'(x)$ had a square in $\alpha_{s_l}$ 
then $W'$ would end in $s_l$ (by the induction hypothesis) and hence $W$  
would not be  an $M$-reduced expression. Hence $W'(x)$ is without squares 
in $\alpha_{s_l}$. 

\vspace{0.2cm}\noindent
Now $W(x)=T_{\alpha_{s_l}}^{Dn_l}(W'(x))$ has a square 
in $\alpha_s$. By Lemma \ref{square} one has: 
\begin{enumerate}
\item either $s_l=s$ and so $W$ ends 
in $s$; 
\item  or else $s_l$ and $s$ are disjoint 
and $W'(x)$ has a square in $\alpha_s$. 
By the induction  hypothesis $W'$ ends in $s$. Since $s_l$ and $s$ commute we switch the 
position of the last two generators  and find that $W$ ends in $s$. 
\end{enumerate}
\end{proof}

\noindent 
For a fixed $a\in A$ the  fundamental group $\pi_1(F(A),p_a^0)$ embeds into 
the groupoid $\pi_1(F(A), S)$. It is also clear that 
$\pi_1(F(A)),p_a^0)$ is kept invariant by the action of an element 
$w\in H(A)$.   

\begin{lem}\label{last}
Assume that the dual intersection graph of $A$ (or, equivalently the surface $F(A)$) 
is connected. 
If $w$ has a nontrivial $M$-reduced expression then $w$ acts non-trivially on 
$O(A)$.    
\end{lem}
\begin{proof}
It is known (see e.g.\cite{CP} and references there) that an $M$-reduced 
expression representing the identity in $H(A)$ is trivial. 
Take then a non-trivial $M$-reduced expression $W$, as above. Since the 
dual intersection  graph of curves is connected there exists some $t\in A$ such 
that  $i(s_l,t)\neq 0$. We will show that $W(\alpha_t)\neq \alpha_t$. 
Since $\alpha_t\in O(A)\subset \pi_1(F(A), p_t^0)$ 
the action of $W$ is nontrivial on $O(A)$. 

\vspace{0.2cm}\noindent
Suppose $W(\alpha_t)=\alpha_t$ and write $W=T_{s_l}^{n_l}W'$. Then 
\[ W'(\alpha_t)=w_{s_l}^{-n_l}(\alpha_t)=T_{s_l}^{-Dn_l}(\alpha_t)\]
Lemma \ref{cal} shows that $T_{s_l}^{-Dn_l}(\alpha_t)$ has a 
square in $\alpha_{s_l}$ and further  from 
Lemma \ref{square} $W'$ ends in $s_l$. 
But then $W$ is not $M$-reduced, contradiction. This proves the claim.  
\end{proof}

\begin{prop}\label{mumu}
Assume that the dual intersection graph of the finite collection 
$A$ is connected, $A$ has at least two elements and $D\geq 2$. 
Then the group  $M(F(A))(A;D)$  
is  a right angled Artin groups of presentation: 
\[ M(F(A))(A;D) 
=\langle T_a^D, a\in A\mid T_a^DT_b^D=T_b^DT_a^D, \mbox {\rm  if } 
a\cap b =\emptyset\rangle \]
\end{prop}
\begin{proof}
Lemma \ref{last} shows that the map $\tau\circ\iota:H(A)\to M(F(A))$ is injective, 
since $M(F(A))$ is a subgroup of the group of automorphisms 
of $\mathbb F$. Therefore $M(F(A))(A;D)$ is isomorphic to $H(A)$, as claimed. 
\end{proof}

\begin{cor}
If $A$ is nontrivial and $\Sigma_{g,d}\setminus F(A)$ has neither disks nor cylinder components  
joining two  distinct boundary  components of $F(A)$ and $D\geq 2$  
then  $M(\Sigma_{g,d})(A,D)$ is a right angled Artin group of presentation: 
\[M(\Sigma_{g,d})(A;D) =\langle T_a^D, a\in A\mid T_a^DT_b^D=T_b^DT_a^D, \mbox {\rm if } 
a\cap b =\emptyset\rangle \]
\end{cor}
\begin{proof}
The embedding $F(A)\subset \Sigma_{g,d}$, with 
$F(A)$ different from a disk or an inessential annulus induces 
a group embedding $M(F(A))\hookrightarrow M(\Sigma_{g,d})$ 
according to (\cite{PR}, Corollary 4.2) 
if and only if $\Sigma_{g,d}\setminus F(A)$ has neither 
disk nor cylinder components joining two  distinct boundary  
components of $F(A)$. Now, since $A$ is nontrivial 
$F(A)$ is neither a disk nor an inessential annulus. 
\end{proof}

{\em End of proof of Theorem \ref{main}.}
It suffices to consider the case when the dual intersection graph of $A$ is 
connected. 
The mapping class group 
$M(\Sigma_{g,d})$  
embeds into ${\rm Aut}(\pi_1(\Sigma_{g,d}, p))$, where 
the base point $p$ is chosen on the boundary $\partial \Sigma_{g,d}$. 
By Lemma \ref{last} for every nontrivial element 
$w\in H(A)$ there is some $z\in O(A)$ such that $\tau(\iota(w))(z)\cdot z^{-1}\neq 1$.  
Since the homomorphism $j:O(A)\to \pi_1(\Sigma_{g,d})$ was assumed 
to be injective it follows that 
$\tau(\iota(w))(j(z))\cdot j(z)^{-1}=j(\iota(w)(z)\cdot z^{-1})\neq 1$. 
Therefore $\tau(\iota(w))$ acts nontrivially on 
$\pi_1(\Sigma_{g,d}, p)$ and thus $\tau(\iota(w))$ is not identity. 
This means that $\tau\circ \iota$ is injective and 
hence the homomorphism  of $H(A)$ onto 
$M(\Sigma_{g,d})(A;D)$ is an isomorphism. 
This settles Theorem \ref{main}.

Let $B=\{c_0,c_1,\ldots,c_{2g}\}$ and 
$C=\{c_1,c_2,\ldots,c_{2g}\}$, where $c_j$ are the curves from the figure 
below. 

\begin{center}
\includegraphics[scale=0.5]{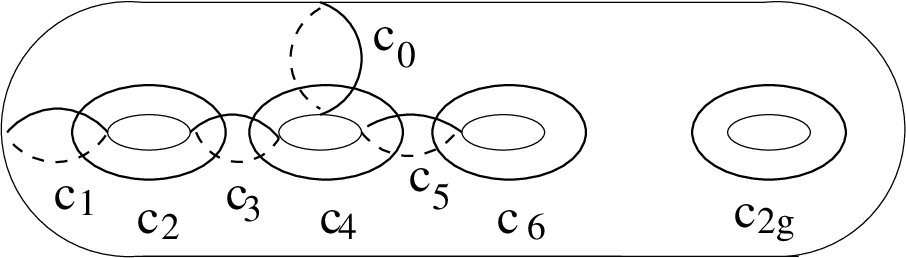}
\end{center}
Let $\Sigma_{g,2}$ and $\Sigma_{g,1}$ be the regular neighborhoods in 
$\Sigma_g$ of the union of curves from $B$ and respectively $C$.  

\begin{cor}
The groups $M(\Sigma_{g,2})(B;D)$  and  $M(\Sigma_{g,1})(C;D)$ 
are right angled Artin groups with the presentations:
\[ M(\Sigma_{g,2})(B;D)
= \langle T_{c_j}^D, j=0,\ldots,2g; T_{c_j}^D T_{c_k}^D=  T_{c_k}^D T_{c_j}^D, 
\mbox{\rm if } j<k, k\neq j+1, (j,k)\neq (0,4)\rangle \]
and respectively:
\[ M(\Sigma_{g,1})(C;D)=\langle T_{c_j}^D, j=1,\ldots,2g; T_{c_j}^D T_{c_k}^D=  T_{c_k}^D T_{c_j}^D, 
\mbox{\rm  if } j<k, k\neq j+1, (j,k)\neq (0,4)\rangle \]
\end{cor}

\begin{proof}   
Here is a direct  simpler proof which uses the proof given in 
\cite{CP} for small Artin groups. 
Let $E_{2g}$ be the Artin group associated to the Dynkin graph of type 
$E_{2g}$, which is the tree whose vertices are in one-to-one correspondence 
with the curves $c_0,c_1,\ldots,c_{2g}$ from the figure above 
and whose edges join two vertices only if the respective curves 
have one intersection point. 
Observe that $A_{2g}$ is the Dynkin subgraph associated to the curves 
$c_1,c_2, \ldots,c_{2g}$.

Let now $E_{2g}[D]$ denote the subgroup of $E_{2g}$ generated 
by  $T_{c_j}^D$, $j=0,1,\ldots,2g$. 
Crisp and Paris proved in \cite{CP} that the subgroup 
$E_{2g}[D]$ has the following  right angled Artin group presentation:  
\[ E_{2g}[D]=\langle T_{c_j}^D, j=0,\ldots,2g \mid 
T_{c_j}^D T_{c_k}^D=  T_{c_k}^D T_{c_j}^D, 
\mbox{\rm  if } j<k, k\neq j+1, (j,k)\neq (0,4)\rangle \]

\vspace{0.2cm}\noindent
The regular neighborhoods  
$F(B)$  and $F(C)$ are homeomorphic to $\Sigma_{g,2}$ and  
$\Sigma_{g,1}$, respectively.

\vspace{0.2cm}\noindent
An essential ingredient of the proof in \cite{CP} is 
the natural representation of the Artin group $E_{2g}$ into the 
mapping class group $M(F(B))$. Consequently $E_{2g}$ acts by automorphisms 
on the fundamental groupoid $\pi_1(F(B); S)$, where 
$S=\{s_0,\ldots,s_{2g}\}$ is a set of boundary base points, 
one base point for each annulus. Set 
$\tau: E_{2g}\to {\rm Aut}(\pi_1(F(B);S))$ for this representation. 

\vspace{0.2cm}\noindent
Let then $H(B)$ and $H(C)$ be the right angled Artin group 
\[ H(B)= \langle a_{j}, j=0,\ldots,2g \mid 
a_{j} a_{k}=  a_{k} a_{j}, 
\mbox{\rm  if } j<k, k\neq j+1, (j,k)\neq (0,4)\rangle \]
\[ H(C)= \langle a_{j}, j=1,\ldots,2g \mid 
a_{j} a_{k}=  a_{k} a_{j}, 
\mbox{\rm if } j<k, k\neq j+1, (j,k)\neq (0,4)\rangle \]
There is a homomorphism $\iota: H(B)\to E_{2g}$ that sends each 
$a_j$ into $T_{c_j}^D$.

The key point of the proof from \cite{CP} is that, given any  
non-trivial element $w\in H(B)$, the automorphism $\tau(\iota(w))$   
acts non-trivially on some element of $\pi_1(F(B),S)$ and 
hence $\tau(\iota(w))\neq 1$. This shows that $\iota$ 
injects $H(B)$  into $E_{2g}$.  

\vspace{0.2cm}\noindent
However this proof also shows that the right angled Artin group 
$H(B)$ injects into the mapping class group $M(F(B))$. The 
corresponding map sends $a_{j }$ into the Dehn twist $T_{c_j}^D$. 
As $M(F(B))$ is  actually $M(\Sigma_{g,2})(B;D)$ the claim follows.

\vspace{0.2cm}\noindent
The same proof works for  the sub-family $C$. 
\end{proof}

\vspace{0.2cm}\noindent
We can slightly generalize the previous results to subgroups 
generated by not necessarily equal powers of Dehn twists. 
\begin{prop}
Let $A$ be a nontrivial sparse collection of curves on $\Sigma_{g,d}$. 
Then the subgroup of $M(F(A))$  generated 
by the powers $T_a^{D(a)}$, where $|D(a)|\geq 2$, $a\in A$ is a right angled 
Artin group.   
\end{prop}
\begin{proof}
The proof from above applies with only minor modifications. 
\end{proof}

\begin{rem}
If $D(a)=D$, for a non-separating curve $a$ and $D(a)=1$, for all other simple closed curves $a$, 
then the subgroup generated by {\em all} the powers $T_a^{D(a)}$  
is the level $D$ subgroup of the mapping class group of $\Sigma_g$, 
namely the kernel of $M(\Sigma_g)\to Sp(2g,\Z/D\Z$). This is proved by 
McCarthy in (\cite{McC}, Theorem 2.8). In particular, in this case the subgroup is 
of finite index. 
\end{rem}

\section{Power subgroups of the mapping class group}

\subsection{ $M_g[D]$ and symplectic groups}\label{firstprop}
We fix once for all a symplectic basis 
$\{a_i,b_i\}_{i=1,\ldots,g}$ in homology consisting of classes 
of simple loops  
and denote by $P:M_g\to Sp(2g,\Z)$ the natural 
homomorphism.

\begin{prop}
If $g\geq 2$ then $P$ sends $M_g[D]$ onto the special congruence subgroup 
\[Sp(2g,\Z)[D]=\ker (Sp(2g,\Z)\to Sp(2g,\Z/D\Z))\]
\end{prop}
\begin{proof}
The action of the Dehn twist $T_b$ in homology is given by 
\[ T_b^D a = a + D \langle a,b\rangle b \]
where $\langle a,b\rangle$ is the algebraic intersection number on $\Sigma_g$. 
Therefore  $T_b^D(a)-a$ belongs to  the submodule 
$D H_1(S_g,\Z)$ of  $H_1(S_g,\Z)$,  
for any $b\in H_1(S_g,\Z)$. 
This implies that $P(T_b^D)\in Sp(2g,\Z)[D])$ and hence 
$P(M_g[D])$ is a normal subgroup of $Sp(2g,\Z)[D]$.

Recall that $Sp(2g,\Z)$ is the group of matrices $A$ with integer 
entries which satisfy $AJA^{T}=J$, where the almost complex 
structure matrix $J$ is the direct sum  of $g$ blocks
$\left(\begin{array}{cc}
0 & 1 \\
-1 & 0 
\end{array}\right)$.

Consider the elementary matrices  
\[ SE_{i\tau(i)}[D]=I_{2g}+ D E_{i\tau(i)}\]
\[ SE_{ij}[D]=I_{2g}+D E_{ij}-(-1)^{i+j}DE_{\tau(j)\tau(i)}\]
where $\tau$ is the permutation $\tau(2j)=2j-1, \tau(2j-1)=2j$, 
for  $1\leq j\leq g$ and $E_{ij}$ denotes the matrix having a single 
non-zero  unit entry at position $(ij)$. 
By direct computation we find that: 
\[ SE_{12}[D]=P(T_{a_1}^{-D})\]
\[ SE_{13}[D]=P(T_{b_2}^{-D}T_{a_1}^{-D}T_c^D)\]
\[ SE_{14}[D]=P(T_{a_2}^DT_{a_1}^DT_{f}^{-D})\]
where $c$ and $f$ are simple closed curves whose homology classes  
are $a_1+b_2$ and $a_1+a_2$ respectively.

Therefore the elementary congruence subgroup of level $D$, which is defined as 
the matrix group generated by the matrices 
$SE_{ij}[D]$, is contained in $P(M_g[D])$. 
Now, a deep result of Mennicke (see \cite{Me1,Me2,BMS}) 
says that the elementary congruence subgroup  
coincides with the congruence subgroup $Sp(2g,\Z)[D]$, 
if $g\geq 2$. Therefore $P(M_g[D])=Sp(2g,\Z)[D]$, as claimed.
\end{proof}

\begin{rem}
If $g=1$ then $M_g[D]$ might be of infinite index in $SL(2,\Z)$ 
(see  \cite{New}). 
\end{rem}

\begin{cor}
The group $M_g[D]$  is torsion-free and consists of pure mapping classes 
when $D\geq 3$ and $g\geq 2$. 
\end{cor}
\begin{proof}
Serre's Lemma tells us that torsion elements in the mapping class group 
act non-trivially on the homology with $\Z/D\Z$ coefficients for any 
$D\geq 3$. 

Recall that a mapping class $h$ is pure 
if $h^n(\gamma)=\gamma$ implies that $h(\gamma)=\gamma$, 
for each isotopy class of a simple closed curve $\gamma$.
Then the second claim is a simple consequence of Ivanov's results (see \cite{I1,I2})
concerning pure classes.
\end{proof}

\subsection{Power subgroups and symplectic groups}\label{powersbgs}

\vspace{0.2cm}\noindent
We start by analyzing the images of the power subgroups in the 
symplectic group. This amounts to find the power subgroups 
of the symplectic group. 
Let $g\geq 2$ and recall that $P$ denotes the  
homomorphism $M_g\to Sp(2g,\Z)$ induced by a homology basis.
We already saw in section \ref{firstprop} that $P(M_g[D])=Sp(2g,\Z)[D]$. 
Moreover since $P$ is surjective $P(X_g[D])$ is a normal subgroup of 
$Sp(2g,\Z)$ containing $Sp(2g,\Z)[D]$.  We have then an obvious 
surjective homomorphism:
\[ L:Sp(2g,\Z/D\Z)=Sp(2g,\Z)/Sp(2g,\Z)[D]\to Sp(2g,\Z)/P(X_g[D])\]
Our first technical result is 
the following:  
\begin{lem}\label{cong}
For any integer $D\not\equiv 0({\rm mod}\,\, 6)$ and any proper ideal $J\subset \Z/D\Z$  
there exists an element in the kernel of $L$ which is 
not central after reduction mod $J$.  
\end{lem}
\begin{proof}
It suffices to find a matrix in $C\in Sp(2g,\Z/D\Z)$ whose 
power $C^D$ is neither the identity ${\mathbf 1}$ nor $-\mathbf 1$ 
modulo the ideal $J$, since the center of $Sp(2g,\Z/D\Z)$ consists of 
$\{{\mathbf 1}, -\mathbf 1\}$ (see \cite{Kl}, Prop.2.1).  
Since $C^D$ belongs to $\ker L$ this will prove the lemma. 
 
\vspace{0.2cm}\noindent
We look for $C$ of the form 
$A\oplus A\oplus \cdots \oplus A$ where 
$A$ is a 2-by-2 matrix. We take a lift of $A$ with integer entries. 
Then $C^D$ has the form 
$A^D\oplus A^D\oplus \cdots \oplus A^D$. 
Since $A\in SL(2,\Z)$ we have 
\[ A^2=t A- {\mathbf 1}\]
where $t$ is the trace of $A$. It follows that  
\[ A^D= Q_{D-1}(t) A - Q_{D-2}(t){\mathbf 1} \]
where $Q_k(t)\in\Z[t]$ are polynomials in the variable $t$ 
determined by the recurrence relation: 
\[ Q_{n}(t)=tQ_{n-1}(t)-Q_{n-2}(t)\]
with initial values $Q_0=1, Q_1(t)=t$. 

\vspace{0.2cm}\noindent
We obtain therefore, by induction on $D$, the following formulas: 
\[ Q_{D-1}(0)=\left\{\begin{array}{rl}
(-1)^{\frac{D-1}{2}}, & \mbox{ if } \,  D\equiv 1({\rm mod}\,\, 2), \\
0, & \mbox{ if } \,  D\equiv 0({\rm mod}\,\, 2).\\
\end{array}
\right.
\] 
\[ Q_{D-1}(-1)=\left\{\begin{array}{rl}
1, & \mbox{ if } \,  D\equiv 1({\rm mod}\,\, 3),  \\
-1, & \mbox{ if } \,  D\equiv 2({\rm mod}\,\, 3), \\
0, & \mbox{ if } \,  D\equiv 0({\rm mod}\,\, 3).\\
\end{array}
\right.
\] 
If the reduction mod $J$ of $C^D$ is trivial for all $C$ 
as above then $Q_{D-1}(t)\equiv 0({\rm mod}\,\, J)$ for all $t$, since  
there exist matrices $A$ of given trace $t$ having some entry 
off-diagonal  which is congruent to $1$ mod $D$, for instance $A=\left(\begin{array}{cc} 
t & 1 \\
-1 & 0\\
\end{array}\right)$.   
Now, either $Q_{D-1}(-1)$ or $Q_{D-1}(0)$ is $\pm 1$ mod $D$, hence 
$J$ is trivial. This proves the claim. 
\end{proof}

\begin{rem}\label{case6} 
The conclusion of Lemma \ref{cong} does not hold when 
$D\equiv 0({\rm mod}\,\, 6)$. For instance 
$Q_{5}(t)=t(t-1)(t+1)(t^2-3)$ and thus $Q_5(t)\equiv 0({\rm mod}\,\, 6)$ 
for every integer $t$. More generally $Q_{6k-1}(t)\equiv 0({\rm mod}\,\, 6)$, 
for every integer $k$. It suffices to observe that: 
\[ Q_{D-1}(1)=\left\{\begin{array}{rl}
1, & \mbox{ if } \,  D\equiv 1({\rm mod}\,\, 6),\mbox{ or } \, D\equiv 2({\rm mod}\,\, 6),  \\
-1, & \mbox{ if } \,  D\equiv 4({\rm mod}\,\, 6),\mbox{ or } \, D\equiv 5({\rm mod}\,\, 6), \\
0, & \mbox{ if } \,  D\equiv 3({\rm mod}\,\, 6),\mbox{ or } \, D\equiv 6({\rm mod}\,\, 6). \\
\end{array}
\right.
\] 
and use the previous computations for $Q_{D-1}(0)$ and $Q_{D-1}(1)$. 
\end{rem}

\begin{rem}
Observe that $Q_n$ is the $n$-th Chebyshev polynomial of the second kind 
\[ Q_n(t)=\frac{\sin (n+1){\rm arcos}(t/2)}{\sin{\rm arcos}(t/2)}\]
which can be given by the explicit formula 
\[ Q_n(t)=\sum_{k=0}^{\left[\frac{n+1}{2}\right]}
(-1)^k\frac{(n-k)!}{k!(n-2k)!}t ^{n-2k}\]
Notice that the usual definition for the Chebyshev polynomial  
uses the variable $x$, where $t=2x$ (see \cite{Rivlin} for more details). 
\end{rem}

\begin{prop}\label{ima}
Suppose that $g\geq 2$, $D$ is of the form $p^m$ 
for a prime $p$, $m\in\Z_+$ and additionally  
$g\geq 3$, $m\geq 2$ when $p\in\{2,3\}$. 
Then $P(X_g[D])$ is all of $Sp(2g,\Z)$. 
\end{prop}
\begin{proof}
We want to prove that the image of $L:Sp(2g,\Z/D\Z)\to Sp(2g,\Z)/P(X_g[D])$,  
(introduced at the beginning of section \ref{powersbgs}) is trivial.  
Since the homomorphism $L$ is surjective, this will prove our claim. 
To this purpose we analyze its kernel $\ker L$.

Now, the normal subgroups of symplectic groups over local rings were 
described by Klingenberg (see \cite{Kl}, Lemma 3.2) and Jehne (\cite{Je}), in the case 
when $D=p^m$, $p$ prime and $p\not\in\{2,3\}$. 
The most general statement can be found in 
(\cite{HO}, Thm. 9.1.7, p.517) where one also considered 
$p\in\{2,3\}$ but $g\geq 3$. The above cited result is that 
under these conditions all normal subgroups of 
$Sp(2g,\Z/D\Z)$ (where $D=p^m$, such that 
$\Z/D\Z$ is a local ring) are {\em congruence} subgroups, 
namely they contain the kernel $Sp(2g, \Z/D\Z)[J]$ of the homomorphism 
$Sp(2g, \Z/D\Z)\to Sp(2g,(\Z/D\Z)/J)$, for some 
ideal $J$. 
This implies that there exists an ideal $J\subset \Z/D\Z$ 
for which $\ker L$ contains  $Sp(2g, \Z/D\Z)[J]$.

On the other hand,  if $J$ were a proper ideal of $\Z/D\Z$ Lemma \ref{cong} would  provide an element  
of $\ker L$ which does not belong to $Sp(2g, \Z/D\Z)[J]$.  
Therefore $J=\Z/D\Z$, whenever  $p\not\in\{2,3\}$ or $m\geq 2$, and hence the map $L$ 
is trivial. 
\end{proof}

\begin{rem}
The projective symplectic group  
$PSp(2g,\Z/DZ)$ is simple when $D$ is prime, except when $g=1, D\in \{2,3\}$ 
(where it coincides with the permutation group $S_3$ and respectively 
the alternating group $A_4$) and $g=2$, $D=2$ (when it coincides with 
the permutation group $S_6$). 
\end{rem}

\begin{rem}
When $g=2$ and $D=2$ the image of $P(X_2[2])$ is of index 2 
in $Sp(4,\Z/2\Z)$. The  subgroup generated by squares of 
elements in $S_6$ is the index 2 alternating subgroup $A_6$. 
In fact any square has even signature and  $A_6$ is also
the commutator subgroup of $S_6$. Observe that $[a,b]= (ab)^2$, if $a^2=b^2=1$  
and commutators of transpositions generate $A_6$. 
Finally we have the exact sequence: 
\[ 1\to \Z/2\Z \to P(X_2[2])\to A_6\to 1.\]
\end{rem}

\vspace{0.2cm}\noindent
In the general case when $D$ is not a power of a prime the image of 
$X_g[D]$ might be  strictly smaller that $Sp(2g,\Z/D\Z)$. This is  
clear when $D\equiv 0({\rm mod}\,\, 6)$, since Remark \ref{case6} shows  
that the image of $P(X_g(D))\subset Sp(2g,\Z/D\Z)$ 
into $Sp(2g,\Z/6\Z)$ must be central. A similar result holds more generally. 
Let us set: 
\[ o_c(D)=\min \{d; A^d\in Z(Sp(2g,\Z/D\Z)), 
\mbox{ for any } A \in Sp(2g,\Z/D\Z)\}\]
where $Z(G)$ stands for the center of the group $G$. 
Write $D$ as $D=q_1 q_2\cdots q_mD'$, where $q_j$ are powers of 
distinct primes and $D'\in\Z$. 
Set $V=\{j; o_c(q_j) \mbox{ divides } D\}\subset \{ 1,2,\ldots, m\}$ and 
$\nu(D)=\prod_{j\in V} q_j$. Consider also the 
general congruence subgroup 
$GSp(2g,\Z/D\Z)[F]$ which is the preimage of $Z(Sp(2g,\Z/F\Z))$
under the reduction mod $F$ homomorphism 
$Sp(2g,\Z/D\Z)\to Sp(2g,\Z/F\Z)$. 
\begin{prop}\label{scg}
The image $P(X_g[D])$ is contained in 
the general congruence subgroup $GSp(2g,\Z)[\nu(D)]$. 
\end{prop}
\begin{proof}
Consider the homomorphism 
$p_j: Sp(2g,\Z/D\Z)\to Sp(2g,\Z/q_j\Z)$ which reduces entries modulo $q_j$. 
If $A\in Sp(2g,\Z/D\Z)$ then $p_j(A^D)$ is central for any 
$A\in Sp(2g,\Z/D\Z)$  if $o_c(q_j)$ divides $D$. 
Therefore the $D$-th power subgroup of 
$Sp(2g,\Z/D\Z)$ is contained into 
$\cap_{j\in V}GSp(2g, \Z/D\Z)[q_j]$, which  
can be identified with $GSp(2g,\Z)[\nu(D)]$. 
\end{proof}

\subsection{Proof of Theorem \ref{nontrivial} }
Theorem \ref{nontrivial} (1) can be restated as follows: 
\begin{prop}\label{nontrivial1}
There exist infinitely many  integers $D$ for which  
$P(X_g[D])$ is a proper subgroup of $Sp(2g,\Z)$, for $g\geq 2$. 
In particular $M_g/X_g[D]$ are non-trivial 
torsion groups, for these values of $D$. 
\end{prop}
\begin{proof}
It is clear that $o_c(q)$ is a divisor of 
the order of $Sp(2g,\Z/q\Z)$, although this upper bound is far 
from being optimal. Let $D={\rm l.c.m.}(o_c(q),q)$. 
Thus we can write $D=qD'$ for some integer $D'$, and we know that  
$o_c(q)$ divides $D$. Therefore $\nu(D)$ is divisible by $q$. Henceforth 
there exist infinitely many  integers $D$ for which  $P(X_g[D])$ is a proper  
subgroup of $Sp(2g,\Z)$, by Proposition \ref{scg}. 
In particular $M_g/X_g[D]$ is a non-trivial 
torsion group. 
\end{proof}

Notice however that $P(X_g[D])$ is always of finite index in $Sp(2g,\Z)$
since it contains the congruence subgroup $P(M_g[D])$. 
The second step in the study of $X_g[D]$ is to understand the interactions 
with the torsion subgroup of $M_g$. We restate here 
Theorem \ref{nontrivial} (2) for the sake of completeness. 

\begin{prop}\label{divisor}
We have $X_g[D]=M_g$, for $g\geq 3$,  
if $4g+2$ does not divide $D$.  
\end{prop}
\begin{proof}
The chain relation (see e.g. \cite{FaMa}, 4.4) shows that whenever $c_1,c_2,\ldots,c_k$ are 
simple closed curves forming a chain i.e. consecutive $c_j$ have a common 
point and are otherwise disjoint, then: 
\begin{enumerate}
\item if $k$ is even we have: 
\[ (T_{c_1}T_{c_2}\cdots T_{c_k})^{2k+2} = T_{d}\]
and also: 
\[ (T_{c_1}^2T_{c_2}\cdots T_{c_k})^{2k}=T_d\]
where $d$ is the boundary of the regular neighborhood of the union of the 
$c_j$. 
\item if $k$ is odd we have: 
\[ (T_{c_1}T_{c_2}\cdots T_{c_k})^{k+1}=T_{d_1}T_{d_2} \]
and respectively: 
\[(T_{c_1}^2T_{c_2}\cdots T_{c_k})^{k}=T_{d_1}T_{d_2} \]
where $d_1,d_2$ are the boundary curves of the regular neighborhood of the union of the 
$c_j$. 
\end{enumerate}

\vspace{0.2cm}\noindent
As a consequence the element $a=T_{c_1}T_{c_2}\cdots T_{c_{2g}}$ is of order 
$4g+2$ and the element  $b=T_{c_1}^2T_{c_2}\cdots T_{c_{2g}}$ is of 
order $4g$, where $c_1,c_2,\ldots,c_{2g}$ are the curves from the first 
figure. 

\begin{lem}\label{normal}
The normal subgroup generated by $a^k$ is $M_g$ when 
$k\leq 2g$ and $g\geq 3$ and of index 2 when $g=2$. 
\end{lem}
\begin{proof} See (\cite{KoHa}, Theorem 4). 
\end{proof} 

\vspace{0.2cm}\noindent
Let $\pi:M_g\to M_g/X_g[D]$ be the projection. 
We have then $a^{4g+2}=1$. Set $k={\rm gcd}(4g+2, D)< 4g+2$. 
In the quotient $M_g/X_g[D]$ we have 
also $\pi(a^{D})=1$ and hence $\pi(a^{k})=1$. We have either 
$k\leq 2g$ or else $k=2g+1$.
 
\vspace{0.2cm}\noindent
If $k < 2g+1$ Lemma \ref{normal} shows that 
the quotient $M_g/X_g[D]$ is trivial. 

\vspace{0.2cm}\noindent
If $k=2g+1$ recall that we have also $b^{4g}=1$  and hence 
$\pi(b)=1$. This implies that 
$\pi(a)=\pi(T_{c_1}T_{c_{2}}\cdots T_{c_g})=\pi(T_{c_1}^{-1})$.

\vspace{0.2cm}\noindent
By recurrence on $k$  we can show that 
$a^k(c_1)=c_{k+1}$, if $k\leq 2g$, where $c_{2g+1}$ is the curve from the figure below: 

\begin{center}
\includegraphics[scale=0.5]{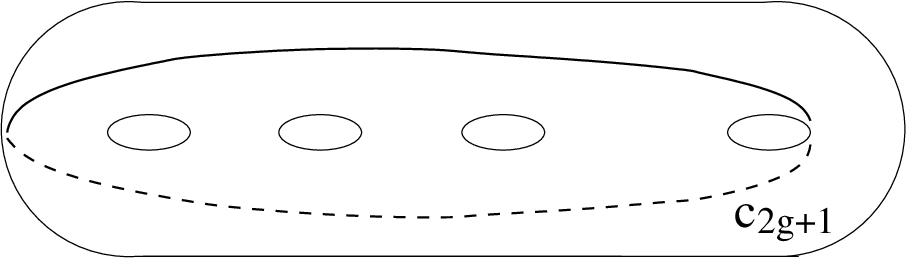}
\end{center}

\vspace{0.2cm}\noindent
Thus  
\[ T_{c_1}^{-1}a^kT_{c_1}a^{-k}=T_{c_1}^{-1}T_{a^k(c_1)}=T_{c_1}^{-1}T_{c_{k+1}}
\]
Therefore 
\[ \pi(T_{c_1}^{-1}T_{c_{k+1}})=\pi(T_{c_1}^{-1}a^kT_{c_1}a^{-k})=1 \]
so that 
\[\pi(T_{c_1})=\pi(T_{c_2})=\cdots =\pi(T_{c_{2g}})\]

\vspace{0.2cm}\noindent
The braid relations in $M_g$ read 
\[ T_{c_0}T_{c_4}T_{c_0}=T_{c_4}T_{c_0}T_{c_4}\]
and 
\[ T_{c_1}T_{c_0}=T_{c_0}T_{c_1}\]
from which one can find 
\[ \pi(T_{c_0})=\pi(T_{c_1})\]
Thus the images by $\pi$ of all  standard $2g+1$ 
generators of $M_g$ coincide and since $H_1(M_g)=1$, for $g\geq 3$, 
we obtain:  
\[ \pi(T_{c_i})=1, \mbox{ for all } i=0,1,\ldots 2g \]
Thus the quotient group is trivial. 
\end{proof}

\begin{rem}
One knows that $M_g/M_g[2]$ is finite (see \cite{Hum}), when $g\geq 2$,   
and $M_g/X_g[2]$ is the further quotient obtained by adjoining 
all squares as relations. 
Thus the quotient is a finite commutative 
2-torsion group. But $M_g$ is perfect (when $g\geq 3$) and  hence
it does not have surjective morphisms into nontrivial abelian groups. Thus  
$M_g/X_g[2]$ should be trivial, for $g\geq 3$.   
\end{rem}

\begin{rem}
For every non-separating curve $d$ we can 
find a chain $c_1,c_2,\ldots, c_{2g-1}$ whose boundary is made of two curves 
isotopic to $d$ and hence 
\[(T_{c_1}^2T_{c_2}\cdots T_{c_{2g-1}})^{2g-1}=T_{d}^2 \]
Since $T_d$ and $T_{c_i}$ commute we have 
\[ ((T_{c_1}^2T_{c_2}\cdots T_{c_{2g-1}})^{1-g}T_d)^{2g-1}=T_d\]
Thus every Dehn twist along a non-separating curve is a $(2g-1)$-power.  
Since  these Dehn twists generate $M_g$ it follows that 
$X_g[2g-1]=M_g$, for $g\geq 3$. 
\end{rem}

\begin{cor}
The index of a normal subgroup of $M_g$ is a multiple of $4g+2$, when $g\geq 3$. 
\end{cor}
\begin{proof}
In fact $X_g[N]$ is contained in a normal subgroup of index $N$.
Proposition \ref{divisor} implies the claim.  
\end{proof}

\subsection{Proof of Theorem \ref{tors}}
For a group $G$ denote by $Q(G)[D]$ the quotient of $G$ by its 
$D$-th power subgroup $X(G)[D]$. 
The key ingredient we shall use is the deep result of 
Adian and Novikov (see  \cite{Adian}), Lys\"enok  (\cite{Ly}) and Sergei Ivanov 
(see \cite{IvanovS}) that the free Burnside group 
$Q({\mathbb F}_2)[D]$ is infinite for large $D$ (e.g. $D\geq 8000$).

\begin{lem}\label{surjective}
If $G\to H$ is surjective then $Q(G)[D]\to Q(H)[D]$ is also surjective. 
\end{lem}
\begin{proof}
It suffices to see that it is well-defined and thus surjective. 
\end{proof}

\begin{lem}\label{subgrou}
If $G\subset H$ is a subgroup of index $n$ and 
$Q(G)[D]$ is infinite then $Q(H)[n!D]$ is infinite.
When $G$ is a normal subgroup then  $Q(H)[nD]$ is infinite. 
\end{lem}
\begin{proof}
If $G$ is a normal subgroup of $H$ then for every $a\in H$ 
we have $a^n\in G$. If $G$ is not necessarily normal 
then we claim that for every $a\in H$ we have $a^{n!}\in G$. 
In fact, by our assumption there are  only $n$ distinct left cosets 
of $G$ in $H$.  Thus the following $(n+1)$ left cosets 
$G$, $aG$, $a^2G, \ldots, a^nG$ cannot be distinct. This means 
that there exists some non-zero integer $p\leq n$ such that $a^p\in G$. 
Since $p$ divides $n!$ it follows that $a^{n!}\in G$, as claimed.  
  
\vspace{0.2cm}\noindent
Therefore $X(H)[nD]\subset X(G)[D]\subset G\subset H$ 
if $G$ is normal and 
$X(H)[n!D]\subset X(G)[D]\subset G\subset H$, otherwise. 
The Lemma follows from this. 
\end{proof}

\begin{lem}\label{gen0}
We have $Q(M_0^n)[n(n-1)(n-2)(n-3)D]$ is infinite if $n\geq 4$ and $D\geq 8000$. 
\end{lem}
\begin{proof}
Observe that $M_0^n$ contains the 
index $n(n-1)(n-2)(n-3)$  subgroup $U$ which preserve point-wise 
four punctures. Let $PM_0^n$ denote the subgroup of pure mapping classes 
in $M_0^n$ which preserve point-wise all punctures. 
Then $U$ surjects onto $PM_0^4$, by forgetting all but the four fixed 
punctures. But  $PM_0^4$ is  isomorphic to 
the free group ${\mathbb F}_2$. 
Thus Lemmas \ref{surjective} and \ref{subgrou} settle the claim. 
\end{proof}

\vspace{0.2cm}\noindent
The proof of Theorem \ref{tors} 
 follows now from the following exact sequence:
\[ 1\to \Z/2\Z \to M_2\to M_0^6\to 1 \]
and Lemmas \ref{surjective} and \ref{gen0}.

\begin{rem}
The same proof shows that the group  $Q(C_{M_g}(j)((2g+2)!D)$ 
associated to the centralizer $C_{M_g}(j)$ of the hyper-elliptic 
involution $j$  is infinite as 
soon as $D$ is large enough.   
\end{rem}

\begin{rem}
One might speculate that for large  values of $D$ the subgroup 
$X_g[g!(4g+2)D]$ is of infinite index in $M_g$ and 
the quotient is a finitely generated torsion group of exponent $g!(4g+2)D$. 
Moreover, in this case it would exist $N(g)$, which divides 
$g!(4g+2)$, such that $Q(M_g)[N(g)D]$ is infinite 
for large enough $D$, while $Q(M_g)[D]$ is finite for  every 
$D$ not divisible by $N(g)$. This would follow if it were  
exist a finite index subgroup of $M_g$ which surjects onto a free 
non-abelian group. 
\end{rem}

{
\small      
      
\bibliographystyle{plain}

}

\end{document}